\documentclass[letterpaper, 10 pt, conference]{ieeeconf}

\usepackage{graphicx}
 \usepackage[pdftex]{color}
\usepackage{amsmath}
\usepackage{amsfonts} 

\usepackage{array}
\newtheorem{teor}{Theorem}[section]
\newtheorem{corr}{Corollary}[section]
\newtheorem{propo}{Proposition}[section]
\newtheorem{lemm}{Lemma}[section]

\newtheorem{rem}{Remark}
\usepackage{soul}
\newcommand{\tp}{^{\top}}

\newcommand{\beq}{\begin{equation}}
\newcommand{\eeq}{\end{equation}}
\newcommand{\bea}{\begin{eqnarray}}
\newcommand{\eea}{\end{eqnarray}}

\newcommand{\bsea}{\begin{subeqnarray}}
\newcommand{\esea}{\end{subeqnarray}}
\newcommand{\nn}{\nonumber}

\newcommand{\qed}{\hfill $\Box$ \vskip 2ex}

\def\bmat{\left[ \begin{array}}
\def\emat{\end{array} \right]}
\DeclareMathOperator{\tr}{tr} 
\DeclareMathOperator{\diag}{diag} 
\DeclareMathOperator{\dd}{diag^{2}} 
\DeclareMathOperator{\ofd}{ofd}
\DeclareMathOperator*{\argmin}{arg\,min}
\definecolor{Royalblue}{cmyk}{1,0.30,0.2,0.2}

\newcommand{\vale}{\color{black}}

\begin{document}

\title{Robust Identification of ``Sparse Plus Low-rank" Graphical Models:\\ An Optimization Approach}

\author{Valentina~Ciccone, Augusto~Ferrante, Mattia~Zorzi \thanks{V. Ciccone, A. Ferrante, and M. Zorzi are with   the Department of Information
Engineering, University of Padova, Padova, Italy; e-mail: {\tt\small valentina.ciccone@dei.unipd.it} (V. Ciccone); {\tt\small augusto@dei.unipd.it} (A. Ferrante); {\tt\small zorzimat@dei.unipd.it} (M. Zorzi).}}

\markboth{DRAFT}{Shell \MakeLowercase{\textit{et al.}}: Bare Demo of IEEEtran.cls for Journals}

\maketitle

\begin{abstract}
Motivated by graphical models, we consider the ``Sparse Plus Low-rank" decomposition of a positive definite concentration matrix- the inverse of the covariance matrix. This is a classical problem for which a rich theory and numerical algorithms have been developed.
It appears, however, that the results rapidly degrade when, as it happens in practice, the covariance matrix must be estimated from the observed data and is therefore affected by a certain degree of uncertainty.
We discuss this problem and propose an alternative optimization approach that appears to be suitable to deal with robustness issues in  the ``Sparse Plus Low-rank" decomposition problem.
The variational analysis of this optimization problem is carried over and discussed.
 \end{abstract}

\section{Introduction}\label{Intro}
Graphical Models are used to provide a graph representation of the relations between random variables.
In particular, Gaussian Graphical models can be used to describe the conditional independence relations between the $m$ components of a zero-mean Gaussian random vector $x$ by mean of 
an \textit{interaction graph} $\mathcal{G}(V_{m}, E_{m})$.
This is an undirected graph where the set $V_{m}$ contains $m$ nodes  and $E_{m}$ is the subset of the pairs of nodes that are directly connected by an edge.
In this representation, the $i$-th node represents the $i$-th component $x_i$ of vector $x$; each edge  represents a dependence relation:  no edge between the nodes $i$ and $j$  indicates that $x_i$ and $x_j$ are conditionally independent, given all the others $x_k$; in more formal terms, for any pair of distinct nodes $i$ and $j$, $(i,j)\notin E_{m}$ if and only if $ x_i \perp x_j| \lbrace x_k \rbrace_{k \neq i,j}$.
This graphical structure is very powerful to represent interdependence of the various components in a complex system with many variables, and for this reason this representation has been used and analyzed in a huge amount of papers in Statistics, Engineering, and Signal processing, to mention only the main applications,\cite{songsiri2010topology,songsiri2010graphical,zorzi2016ar,ARMA_GRAPH_AVVENTI,maanan2017conditional,ID_DAHLHAUS,brillinger1996remarks,7402835,Chandrasekaran_latentvariable}.

 The case of graphs with a small number of edges is very interesting in applications because such a pattern highlights a simple structure where the values of most pairs of variables does not have a {\em direct} influence on each other.
Moreover, this structure clarifies the paths of interdependence along which the value of each variable may affect the value of each other.
In most situations, however, the graph is complete (or almost complete) as there is no pairs of variables that are  conditionally independent given the others.
Of course, this {\em may} reveal that we are considering a genuinely complex system where each variable has direct influence on each other, but it may also point to a different and much more interesting situation.
The essence of this situation is well described by the following very simple example.
Let 
$$x=\begin{bmatrix} x_1\\x_2\\\vdots\\ x_m\end{bmatrix}=\begin{bmatrix} y_1+y_{m+1}\\ y_2+y_{m+1}\\\vdots\\ y_m+y_{m+1}\end{bmatrix},$$ where $y_i$, $i=1,\dots,m+1$ are independent Gaussian random variables.
It is immediate to see that the graph associated with $x$ is complete (i.e. each pair of nodes is connected by an edge). 
The interdependence pattern between the variables of the vector $x$, however, has a very interesting and simple structure providing a powerful interpretation.
This structure can be highlighted by considering the {\em augmented}
vector $\tilde{x}:=[\, x^T\; (y_{m+1})^T\,]^T$; 
indeed, it is easy to see that the  graph associated with $\tilde{x}$ has only the $m$ edges connecting $y_{m+1}$ to each of the $x_i$ which provides the following  interpretation: the interdependence between all the \textit{manifest} (observed) variables $x_i$ is completely explained by the dependence of each $x_i$ with a common \textit{latent} (hidden) variable $y_{m+1}$. 
It is now clear why a great effort has been dedicated to uncover this hidden structure by only observing the manifest variables in $x$.

\begin{figure}[htb]
	\centering
	\includegraphics[clip,width=0.4\textwidth]{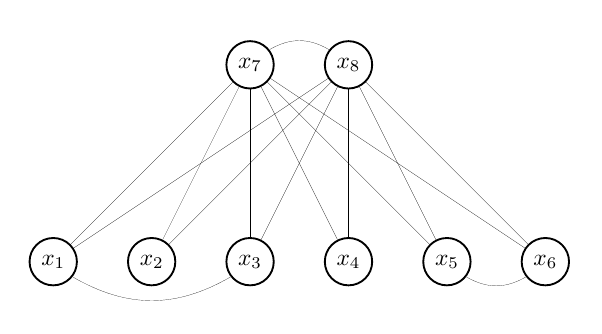}
	\caption{Example of a latent-variable graphical model: the nodes $x_7, \, x_8$ represent the latent variables while the nodes $x_1, \, x_2, \,...,x_6 $ represent the manifest variables.}
\label{fig:figure_two}
\end{figure}

In this paper, we are interested in Gaussian Graphical models with a two-layer structure: our aim is to integrate the set of manifest variables | that we think 
to be arranged in a bottom layer | with a {\em small} number of latent variables 
| arranged in a top layer | in such a way that a drastic reduction in the number of edges between the  manifest  variables is achieved.
An example of this structure is depicted in Figure \ref{fig:figure_two}
where it is easy to see that if we considered only the manifest variables $x_1,\dots,x_6$ the corresponding graph would be complete. However, almost all of the interdependence between these variables is explained by the two latent variables $x_7, \, x_8$ so that an illuminant structure emerges when we integrate these two variables with the observed ones. It is worth noting that this idea has been exploited also for Bayesian networks, \cite{zorzi2017sparse}.

As an application of the theory developed in Dempster seminal paper \cite{Dempster-72},
an identification procedure for such  Gaussian Graphical models has been developed which is based on the  \textit{``Sparse plus Low-rank''} decomposition of the manifest concentration matrix $\Sigma_m^{-1}$ of the manifest vector $x$   ($\Sigma_m$ being the covariance matrix of $x$), \cite{Chandrasekaran_latentvariable,zorzi2016ar}:
\beq\label{splrd}
\Sigma_m^{-1}=S-L.
\eeq
Indeed, such a decomposition, where $S$ is symmetric and positive definite and $L$
is symmetric and positive semi-definite, provides a two-layers graphical model with ${\rm rank}[L]$ latent variables and a number of edges between the observed variables that is equal to the number
of non-zero off-diagonal entries of the matrix $S$.
Therefore, for the solution of our problem we seek for a decomposition of the form \eqref{splrd}  where the rank of the matrix $L$ and the number of non-zero entries of the matrix $S$ are minimized. 

\begin{figure}[thb]
	\centering
	\hspace{-3mm}\includegraphics[width=0.5\textwidth]{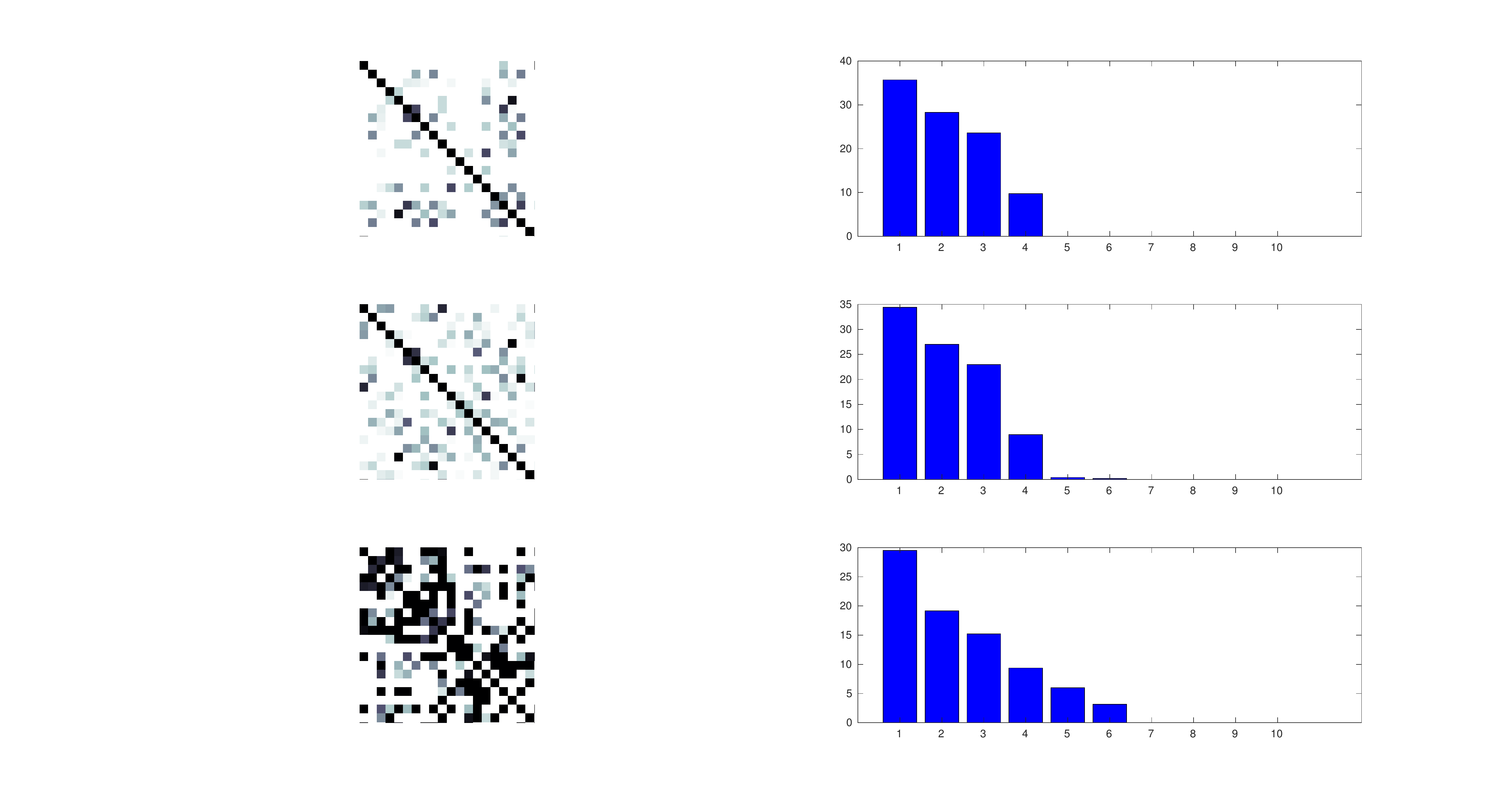}
	\caption{First row: sparsity pattern of the true sparse matrices $S_0$ (and eigenvalues of the true low rank matrix $L_0$. Second and third row: sparsity pattern of the matrices $\tilde{S}_0$ and $\hat{S}_0$ and first $10$ eigenvalues of the  matrices $\tilde{L}_0$ and $\hat{L}_0$ obtained by applying the standard decomposition algorithm to the true covariance matrix $\Sigma_m$ and to the estimated covariance matrix $\hat{\Sigma}_m$, respectively.}
\label{fig:figure_one}
\end{figure}

This classical approach is, however, based on a very artificial assumption: the covariance matrix  $\Sigma_m$ of $x$ is assumed to be known.
On the contrary $\Sigma_m$ is normally estimated from the data so that only a noisy version $\hat{\Sigma}_m$ of $\Sigma_m$ is usually available.
On the other hand,    
the accuracy of this estimation may severely affect the goodness of the result | in terms of minimum rank and maximum sparsity | of the aforementioned optimization problem.
More precisely, even in the case where the data are indeed produced by a mechanism in which a few non-observable variables explain most of the interdependence between the observed variables, relatively small variations of the covariance matrix $\hat{\Sigma}_m$ from the true value $\Sigma_m$ may produce significant changes in the numerical rank of $L$ and the numerical sparsity of $S$. 
To see this, we considered
a sparse matrix $S_0$ of dimension $20$ and a positive semidefinite matrix $L_0$ of dimension $20$ and rank $4$ such that 
$S_0-L_0\succ 0$;
then we considered the matrix $\Sigma_m:=(S_0-L_0)^{-1}$.
We generated a sample of $N=1000$ independent realizations of a Gaussian random vector with  zero mean and covariance matrix $\Sigma_m$ and from them we estimated the sample covariance $\hat{\Sigma}_m$. Then we computed the ``Sparse plus Low-rank'' decompositions  $\Sigma_m^{-1}=\tilde{L}_0-\tilde{S}_0$ and $\hat{\Sigma}_m^{-1}=\hat{L}_0-\hat{S}_0$, using 
the available algorithm in the literature \cite{Chandrasekaran_latentvariable,zorzi2016ar}.
In Figure \ref{fig:figure_one} the sparsity pattern of $S_0$, $\tilde{S}$ and $\hat{S}$ and the first 10 eigenvalues of $ L_0, \tilde{L}_0$ and $\hat{L}_0$ are depicted providing evidence of the degradation of the solution when $\hat{\Sigma}_m$ is substituted to $\Sigma_m$. In fact, it is apparent that when the true covariance matrix is employed the algorithm recover a solution with the correct numbers of latent variables and  of non-zero elements of the matrix  $S$, while when the covariance is estimated (even from as many as $N=1000$ data) the ``sparse plus low-rank'' structure is completely lost.
Therefore, we are dealing with a very delicate problem where the solution is highly sensible to the observed data:  in this work we propose a method to deal with this problem by taking the uncertainty in the estimation of $\Sigma_m$ into account
 as it has been done in robust state estimation paradigms proposed in \cite{levy2004robust,ZORZI2017133}. This leads to a more complex optimization problem: given $\hat{\Sigma}_m$, we propose to {\em compute} the matrix $\Sigma_m$ in such a way that in the decomposition $\Sigma_m^{-1}=S-L$ the rank of $L$ is minimized while the sparsity of $S$ is maximized under a constraint which limits the Kullback-Leibler divergence
between $\Sigma_m$ and the sample covariance  $\hat{\Sigma}_m$ to a prescribed tolerance depending on the precision of $\hat{\Sigma}_m$. {\vale It is worth noting that a similar problem has been addressed with different approaches in \cite{tao2011recovering}, \cite{zhou2011godec}.}

The paper is organized as follows. In Section \ref{sec_pb_form} Gaussian Graphical models and the \textit{Sparse plus Low-rank} structure for the inverse covariance matrix are introduced. Then, the main optimization problem is stated.  
In Section \ref{delta} the choice of the upper bound on the divergence between $\Sigma_m$ and $\hat{\Sigma}_m$ is discussed. Then, in Section \ref{dual_section} we define the dual problem and in Section \ref{existence_uniq} we establish and existence and uniqueness result for the solution of such a dual problem and we show how to obtain the solution of the primal problem.

{\em Notation:} Given a vector space $\cal V$ and a subspace ${\cal W}\subset {\cal V}$, we denote by
${\cal W}^\bot$ the orthogonal complement of ${\cal W}$ in ${\cal V}$.
Given a matrix $M$, we denote its transpose by $M\tp$; if $M$ is a square matrix 
$\tr(M)$ denotes its trace;
moreover, $|M|$ denotes the determinant of $M$ and $\sigma(M)$ denotes the spectrum of $M$, that is, the set of its eigenvalues. 
We endow the space of square real matrices with the following inner product: for $A,B\in{\mathbb R}^{n\times n}$, $\langle A,B\rangle:= \tr(A^TB)$. The Frobenius norm of $M$ is denoted by $\Vert M \Vert$.
The symbol $\mathbf{Q}_n$ denotes the vector space of real symmetric matrices of size $n$. If $X\in\mathbf{Q}_n$ is positive definite or positive semi-definite we write $X\succ 0$ or $X\succeq 0$, respectively.
Moreover, we denote by $\mathbf{D}_n$ the vector space of  diagonal matrices of size $n$; $\mathbf{D}_n$ is clearly a subspace of $\mathbf{Q}_n$. $\mathbf{M}_n:=\mathbf{D}_n^\perp$ is the orthogonal complement of $\mathbf{D}_n$ in $\mathbf{Q}_n$ with respect to the inner product just defined. We denote by $\diag(\cdot)$ both the operator mapping $n$ real elements $d_i, i=1,...,n$ into the diagonal matrix having the $d_i$'s as elements in its main diagonal and 
the operator mapping a matrix $M\in\mathbb{R}^{n\times n}$ into an $n$-dimensional vector containing the diagonal elements of $M$. Then $\diag\diag(\cdot)$, that we denote by
 $\dd(\cdot)$, is the (orthogonal projection) operator mapping  a square matrix $M$ into a diagonal matrix of the same size having the same main diagonal of $M$.
We denote by $\ofd(\cdot)$ the self-adjoint operator orthogonally projecting $\mathbf{Q}_n$ onto $\mathbf{M}_n$, i.e. if $M\in  \mathbf{Q}_n$, $\ofd(M)$ is the matrix of $\mathbf{M}_n$ in which each off-diagonal element is equal to the corresponding element of $M$ (and each diagonal element is zero). We denote by $\otimes$ the Kronecker product between two matrices   and by vec(X) the vectorization of a matrix X formed by stacking the columns of X into a single column vector. Finally, the cardinality of a set $\mathcal{S}$ is denoted by $|\mathcal{S}|$.

\section{Problem Formulation}\label{sec_pb_form}

\subsection{Latent-Variable Gaussian Graphical Models}

Let $\tilde{x}$ be a zero-mean Gaussian random vector 
of dimension $m+l$, that is $\tilde{x}:=[x^\top y^\top]^\top$, where $x:=[x_1 \quad ... \quad x_m]^\top$  plays the role of the manifest vector and $y:=[x_{m+1} \quad ... \quad x_{m+l}]^\top$ plays the role of the latent vector.
We partition the covariance $\Sigma\in \mathbf{Q}_{m+l}$ of $\tilde{x}$ conformably with the partition of $\tilde{x}$ as:
\begin{equation}
 {\Sigma}= 
\begin{bmatrix}
      \Sigma_m  & \Sigma_{lm}^\top  \\
       \Sigma_{lm} & \Sigma_l \\
\end{bmatrix}.
\end{equation}

We are interested in the Gaussian Graphical model of $\tilde x$ and hence in the conditional independence {\vale relations} in $\tilde x$.
Thus we recall a fundamental result stating that two distinct elements of a Gaussian random vector are conditionally independent given all the others if and only if the corresponding element in the concentration matrix (the inverse of the covariance) is zero {\vale \cite{Dempster-72}}. In our case, this reads:
$$
\forall i\neq j, \quad x_i \perp x_j| \lbrace x_k \rbrace_{k \neq i,j} \Leftrightarrow [\Sigma^{-1}]_{ij}=0.
$$

Let us denote by $K$ the inverse of covariance matrix $\Sigma$. Then, $K$ can be also partitioned  as: 
\begin{equation}
K=\Sigma^{-1} = 
\begin{bmatrix}
      K_m  & K_{lm}^\top  \\
       K_{lm} & K_l \\
\end{bmatrix}
\end{equation}
and, by using the \textit{Schur complement}, we can obtain the relationship
\begin{equation}
\label{Schur}
\Sigma_m^{-1}=K_m-K^\top_{lm} K_l^{-1}K_{lm}
\end{equation}
where the sparsity pattern of $K_m$ provides the relations of conditional independence between the manifest variables in $x$ and the rank of $K^\top_{lm} K_l^{-1}K_{lm}$ provides (an upper bound for) the number $l$ of latent variables.

Now we recall that our only data are those of the vector $x$ and hence we may estimate
$\Sigma_m$ while the rest of $\Sigma$ is a purely artificial construction. The previous argument, however,  provides a procedure for identifying a model $\Sigma$ from $\Sigma_m$.
In fact, if we decompose $\Sigma_m^{-1}$ as
\[\Sigma_m^{-1}= S-L, \qquad S\succ 0,\ L\succeq 0,\]
where the rank $l$ of $L$ is as small as possible and $S$ is as sparse as possible
(of course there is a trade off between these two conditions) we may
identify $S$ with $K_m$ and $L$ with $K^\top_{lm} K_l^{-1}K_{lm}$. 

This argument naturally leads to an optimization problem where the optimal $S$ and $L$ must minimize a combination of 
a distance or divergence function $d(\cdot,\cdot)$ between $\Sigma_m$ and $S-L$ and two penalty functions inducing sparsity and low-rankness on $S$ and $L$, respectively, \cite{Chandrasekaran_latentvariable,zorzi2016ar}: 

\begin{equation}
\label{classic}
\begin{aligned}
(\tilde{S}, \tilde{L})=  \argmin_{S,L\in \mathbf{Q_m}} \quad &d(\Sigma_m,S-L)+ \lambda(\gamma \phi_1(S)+\phi_* (L)) \\
\text{subject to } \quad & S-L {\vale \succ  0, \; L\succeq 0}
\end{aligned}
\end{equation}
where $\lambda>0$ and the regularizer is a combination of two penalties: $\phi_1$  which induces sparsity on (the off-diagonal of) $S$ and $\phi_*$ which induces low-rankness on $L$; while $\gamma>0$ plays the role of the balance term. As in \cite{Chandrasekaran_latentvariable,zorzi2016ar}, a natural choice for the regularizers is the following.
\begin{itemize}
\item The \textit{sparsity regularizer} for a matrix $Y\in\mathbf{Q}_n$ is given by an $l_1$-like function:
\[\phi_1(Y)= h_1(Y)=\sum_{k>h}|Y_{hk}|.\]
\item The \textit{low-rank regularizer} is the nuclear norm, which for positive semidefinite matrices can be surrogated by the trace, so that we set:
\[\phi_*(Y)=\tr(Y).\]
\end{itemize}

\subsection{Robust Sparse plus Low Rank Identification}

In practical applications, however, the matrix $\Sigma_m$ is unknown and needs to be estimated from the observed data that we assume to be $N$ independent realizations $x_{(i)}$ of $x$. The typical choice is the sample covariance matrix 
{\vale $\hat{\Sigma}_m=N^{-1} \sum_{i=1}^{N}x_{(i)} x_{(i)}^\top.$ }
As discussed in the Introduction, however, when replacing $\Sigma_m$ with $\hat{\Sigma}_m$ in \eqref{classic} the corresponding solution may rapidly degrade.
To deal with a similar problem in the context of Factor Analysis we have proposed in \cite{NF_CDC,ciccone2017factor} 
an optimization technique that may be adapted to the present setting as discussed below.

Let  $\hat{\Sigma}_m$ be given. We assume  that the ``actual'' $\Sigma_m$ belongs to a ball centred in $\hat{\Sigma}_m$:
\[\mathcal{B}:=\lbrace \Sigma_m \in\mathbf{Q}_m \; s.t. \; \Sigma_m\succ 0, \mathcal{D}_{KL}(\hat{\Sigma}_m \Vert \Sigma_m)\leq \delta/2 \rbrace\]
where $\delta/2$ is the prescribed tolerance and $\mathcal{D}_{KL}(\hat{\Sigma}_m \Vert \Sigma_m)$ is the Kullback-Leibler divergence defined as:
\begin{equation}
\label{KL}
\mathcal{D}_{KL}(\hat{\Sigma}_m||\Sigma_m):= \dfrac{\log |\Sigma_m|-\log |\hat{\Sigma}_m|+\tr(\Sigma^{-1}_m\hat{\Sigma}_m)-m}{2}.
\end{equation}

Therefore, we consider the following problem:
\begin{equation}
\label{problem0}
\begin{aligned}
(\tilde{S}, \tilde{L})=  \argmin_{S,L\in \mathbf{Q_m}} \quad & \tr(L)+\gamma h_1(S)\\
\text{subject to } \quad 
						 & L\succeq 0 \\
						 & S-L \succ 0  \\
						 & \Sigma_m^{-1}= S-L \\
						 &  2 \mathcal{D}_{KL}(\hat{\Sigma}_m||\Sigma_m)\leq \delta.
\end{aligned}
\end{equation}

To streamline the notation let us denote by $X$ the inverse of the matrix $\Sigma_m$, that is $X:=\Sigma_m^{-1}= S-L$.\\ 
Then, the minimization problem in \eqref{problem0} is equivalent to:
\begin{equation}
\label{problem}
\begin{aligned}
(\tilde{S}, \tilde{X})=  \argmin_{S,X\in \mathbf{Q_m}} \quad & \tr(S-X)+\gamma h_1(S)\\
\text{subject to } \quad 
						 & S-X\succeq 0 \\
						 & X \succ 0\\
						 & 2 \mathcal{D}_{KL}(\hat{\Sigma}_m||X^{-1})\leq \delta
\end{aligned}
\end{equation}
where $2 \mathcal{D}_{KL}(\hat{\Sigma}_m||X^{-1})$ is given by:
\[2\mathcal{D}_{KL}(\hat{\Sigma}_m|| X^{-1}):=-\log |X| -\log |\hat{\Sigma}_m|+\tr(X\hat{\Sigma}_m)-m.\]

\subsection{Negative log-likelihood approach}
Before entering into the study of problem \eqref{problem}, it is worth observing that there is a different route leading to essentially the same problem; in fact, it provides an interesting alternative interpretation of our proposed approach.

Let $X_N:=(x_{(1)}, ..., x_{(N)})$ be an i.i.d. sample drawn from $\mathcal{P}(\Sigma_m):=\mathcal{N}(0,\Sigma_m)$. The log-likelihood function is
\begin{align*}
\log p(X_N|\Sigma_m)= & - \dfrac{Nm}{2}\log(2\pi)-\dfrac{N}{2}\log |\Sigma_m| \\
& -\dfrac{1}{2}\sum_{k=1}^{N}x_k^\top \Sigma^{-1}_m x_k.
\end{align*}

By using the so-called ``trace-trick'' we can rewrite it as:
$
\log p(X_N|\Sigma_m)= {\vale -} \dfrac{N}{2}[ m\log(2\pi)+ \log |\Sigma_m|  + \tr(\hat{\Sigma}_m\Sigma^{-1}_m)].
$
Thus, the negative log-likelihood is:
\begin{equation}
\label{likelihood}
l(X_N|\Sigma_m)=  \dfrac{N}{2}[\log |\Sigma_m|+\tr(\hat{\Sigma}_m\Sigma^{-1}_m) + m\log(2\pi)].
\end{equation}
The Kullback-Leibler divergence \eqref{KL} and the negative log-likelihood \eqref{likelihood} differ only from a term not depending on $\Sigma_m$. Accordingly, imposing an upper bound $\delta$ on the Kullback-Leibler divergence is equivalent to impose an upper bound $\bar{l}$ on the desired negative log-likelihood. This leads to the equivalent problem:
\begin{equation*}
\label{problem_likelihood}
\begin{aligned}
(\tilde{S}, \tilde{X})=  \argmin_{S,X\in \mathbf{Q}_m} \quad & \tr(S-X)+\gamma h_1(S)\\
\text{subject to } \quad 
						 & S-X\succeq 0 \\
						 & X \succ 0\\
						 & l(X_N|\Sigma_m) \leq \bar{l}
\end{aligned}
\end{equation*}
where $ \bar{l}:= \dfrac{N}{2}(\delta+\log|\hat{\Sigma}|+m +m\log(2\pi)) $ .


\section{The Choice of $\delta$}\label{delta}
The first issue we have to address in the study of problem \eqref{problem} is how to determine the allowed tolerance $\delta$ that must be selected by taking into account the accuracy of the estimate $\hat \Sigma_m$  which, in turn, depends on the numerosity of the sample size. This can be achieved by choosing a probability $\alpha\in (0,1)$ and a neighborhood of ``radius'' $\delta_\alpha$ (in the Kullback-Leibler topology) centered in $\hat \Sigma_m$ which contains the ``true'' $\Sigma_m$ with probability $\alpha$. An effective approach, based on standard Monte Carlo method, to estimate such $\delta_\alpha$ has been proposed in \cite{ciccone2017factor}.

If the level of $\alpha$ that has been chosen is too large with respect to the numerosity of the available data $N$, the computed $\delta_\alpha$ may turn out to be excessively large so that there exist diagonal matrices $\Sigma_m$ such that $2\mathcal{D}(\Sigma_m || \hat{\Sigma}_m) \leq \delta_\alpha$. In this case problem \eqref{problem} admits the trivial solution $L=0$ and $S$ diagonal. From now on we assume that we are not in this trivial situation or, equivalently, that  $\delta$ in \eqref{problem} is strictly smaller than a certain upper bound $\delta_{max}$ that can be computed by solving the minimization problem:

\begin{align}
\label{min_diagonal}
\delta_{max}:=\min_{D\in\mathbf{D}_n, \; D\succ 0 }  2\mathcal{D}_{KL}(\hat{\Sigma}_m \Vert D ).
\end{align}
The next result provides the solution to this problem.
\begin{propo}
\label{prop_delta_max}
The optimal $D$ solving problem \eqref{min_diagonal} is given by
$D^{opt}= \text{diag}^2(\hat{\Sigma}_m)$,
so that
\begin{equation}
\label{delta_MAX}
\delta_{max}= 2\mathcal{D}_{KL}(\hat{\Sigma}_m \Vert D^{opt})= \log|\hat{\Sigma}_m^{-1} \dd(\hat{\Sigma}_m)|.
\end{equation} 
\end{propo}
{\vale The proof is an easy computation and it is left to the reader.}

\section{The Dual Problem} \label{dual_section}

We reformule the constrained optimization problem in \eqref{problem} as an unconstrained minimization problem using the duality theory.
The Lagrangian of \eqref{problem} is 
\begin{equation}
\label{langrangian}
\begin{aligned}
\mathcal{L}(X,S, U &,\lambda ) = \tr(S)-\tr(X)+\gamma h_1(S)-\tr(U(S-X))  \\
& +\lambda  \left( -\log |\hat{\Sigma}_m|-\log |X|+\tr(X\hat{\Sigma}_m)-m-\delta
 \right) \\
& = \langle S, I-U \rangle +\gamma h_1(S) + \langle X, U-I +\lambda  \hat{\Sigma}_m \rangle  \\
& +\lambda  \left( -\log |\hat{\Sigma}_m|-\log |X|-m-\delta
 \right)  \\
\end{aligned}
\end{equation}
where $\lambda \in \mathbb{R}, \; \lambda \geq 0$, and $U \in \mathbf{Q}_m, \; U\succeq 0$ are the Lagrange multipliers. Then, the dual function is the infimum of $\mathcal{L}$ over $X$ and $S$.\\ 
$\bullet$ \textit{Partial minimization over S.}\\ $\mathcal{L}$ depends on $S$ solely through the terms
\begin{equation}
\label{first_partial_min}
\gamma h_1(S)- \langle S, U-I \rangle.
\end{equation}
As in \cite{songsiri2010topology}, the nonlinear term does not depend on the main diagonal elements of $S$. Therefore, the minimization over the diagonal elements is unbounded below unless $\dd (U-I)=0$, i.e.
\begin{equation}
\label{diag_u}
  \quad U_{ii}=1, \; i=1,...,m. 
\end{equation}
The minimization over the off-diagonal entries of $S$ translates into an independent minimization of

\[ -\left( U_{ij}S_{ij}+S_{ji}U_{ji} \right) + \gamma |S_{ij}| \]
 for each element $i,j$ with $j>i$, which is unbounded below unless:
\begin{equation}
\label{ofd_u}
|U_{ij}|\leq \dfrac{\gamma}{2}, \quad i\neq j.
\end{equation} 
If \eqref{ofd_u}, \eqref{diag_u} hold the infimum of \eqref{first_partial_min}
is equal zero. Therefore, the partial minimization of the Lagrangian over $S$ is
\[ \small
\inf_S \mathcal{L}=\hspace{-0.1cm} 
\begin{cases}
\langle X, \ofd(U) +\lambda \hat{\Sigma}_m \rangle - \lambda(\log |\hat{\Sigma}_m|+\log |X| +m+\delta)\\ \hfill \text{if } \eqref{diag_u}, \eqref{ofd_u} \text{ hold};\\
-\infty \hfill \text{otherwise}
\end{cases} 
\]
$\bullet$ \textit{Partial minimization over X.}\\
If \eqref{diag_u} and \eqref{ofd_u} hold, $\mathcal{L}$ depends on $X$ only through 
\[\langle X, \ofd(U)+\lambda \hat{\Sigma}_m\rangle - \lambda \log |X|\]
which is bounded below if and only if
\begin{equation}
\label{inf_k_1}
\ofd(U) + \lambda \hat{\Sigma}_m \succ 0
\end{equation}
which implies
\begin{equation}
\label{inf_k_2}
\lambda >0.
\end{equation}
If \eqref{inf_k_1} holds, by taking convexity into account, the matrix $X$ achieving the minimum is easily obtained by annihilating the first derivative which yields
\begin{equation}
\label{opt_X}
X=(\lambda^{-1}\ofd(U)+\hat{\Sigma}_m)^{-1}
\end{equation}
and the minimum is therefore  given by
\[\lambda m +\lambda \log |\lambda^{-1}\ofd(U)+\hat{\Sigma}_m|.\]

Therefore, the result of the minimization of the Lagrangian over $S$ and $X$ is:
\begin{equation} \small
\inf_{S,X} \mathcal{L}= 
\begin{cases}
\lambda \log |\lambda^{-1}\ofd(U)+\hat{\Sigma}_m|+ \lambda( -\log |\hat{\Sigma}_m|-\delta)\\ 
 & \hspace{-3.0cm}\text{if } \eqref{diag_u}, \eqref{ofd_u}, \eqref{inf_k_1}, \eqref{inf_k_2} \hbox{ hold}; \\
-\infty & \hspace{-0.7cm}\text{otherwise.}
\end{cases}\nn 
\end{equation}

Let us define the convex, closed and bounded set $\mathcal{U}$ as:
\[\mathcal{U}:=\lbrace U: U=U^\top\succeq 0, \dd (U)=I, |U_{hk}|\leq\dfrac{\gamma}{2}, k\neq h \rbrace.\]
Then, the dual problem is
\begin{equation}
\label{dual_max}
\max_{(U,\lambda) \in \mathcal{C}} \quad \lambda \log |\lambda^{-1}\ofd(U)+\hat{\Sigma}_m| - \lambda(\log |\hat{\Sigma}_m|+\delta)
\end{equation}
where the set $\mathcal{C}$ is defined as:
\[\mathcal{C}:=\lbrace (\lambda, U): U\in\mathcal{U}, \lambda>0, (\lambda^{-1}\ofd(U)+\hat{\Sigma}_m)\succ 0 \rbrace.\]

\section{Existence and Uniqueness} \label{existence_uniq}
We now address existence and uniqueness of the optimal solution of \eqref{dual_max}.
We first reformulate the dual problem as the equivalent minimization problem:
\begin{equation}
\label{dual}
\min_{(U,\lambda) \in \mathcal{C}} \quad \tilde{J}
\end{equation}
with $\tilde{J}$ being the opposite of the dual functional:
\[\tilde{J}:= - \lambda ( \log |\lambda^{-1}\ofd(U)+\hat{\Sigma}_m|  -\log |\hat{\Sigma}_m|-\delta).\]

\subsection{Existence}
To show that \eqref{dual} admits a solution we are going to show that we can restrict our set $\mathcal{C}$ to a compact set $\mathcal{C}_F\subset \mathcal{C}$ over which the minimization problem is equivalent.\\

\begin{lemm} \label{lemma_l0}
Let $(\lambda_k, U_k)_{k\in \mathbb{N}}$ be a sequence of elements in $\mathcal{C}$ such that
\[\lim_{k\rightarrow \infty} \lambda_k= 0.\]
Then $(\lambda_k, U_k)_{k\in \mathbb{N}}$ is not an infimizing sequence for $\tilde{J}$.
\end{lemm}

{\vale The proof is omitted for brevity, see Lemma 5.1 in \cite{ciccone2017factor} for a similar result.}

As a consequence, minimizing the dual functional over the set $\mathcal{C}$ is equivalent to minimize over the set:
\[\mathcal{C}_1:=\lbrace (\lambda, U): U\in\mathcal{U}, \lambda \geq\epsilon, (\lambda^{-1}\ofd(U)+\hat{\Sigma}_m)\succ 0 \rbrace\]
for a certain $\epsilon>0$.\\

The next result allows to further restrict $\mathcal{C}_1$ to a set where  $\lambda$ is bounded.
\begin{lemm}
Let  $(\lambda_k, U_k)_{k\in \mathbb{N}}$ be a sequence of elements in $\mathcal{C}_1$ such that 
\[\lim_{k\rightarrow \infty} \lambda_k= \infty.\]
Then $(\lambda_k, U_k)_{k\in \mathbb{N}}$ cannot be an infimizing sequence for $\tilde{J}$.
\end{lemm}

{\vale The proof is omitted for brevity, see Lemma 5.2 in \cite{ciccone2017factor} for a similar result.}

As a consequence, we can further restrict set set $\mathcal{C}_1$ to the set:
\[\mathcal{C}_2:=\lbrace (\lambda, U): U\in\mathcal{U}, \xi \geq\lambda \geq\epsilon , (\lambda_k^{-1}\ofd(U_k)+\hat{\Sigma}_m)\succ 0 \rbrace\]
for a certain $\xi>0$.\\

Finally, we consider a sequence $(\lambda_k, U_k)_{k\in \mathbb{N}}$  such that as $k\rightarrow \infty$ the minimum eigenvalue of $\lambda^{-1}\ofd(U)+\hat{\Sigma}_m$ tends to zero. This implies $|\lambda_k^{-1}\ofd(U_k)+\hat{\Sigma}_m|\rightarrow 0$ and hence $\tilde{J}\rightarrow + \infty$. Therefore, such a sequence cannot be an infimizing sequence.

In conclusion, we can restrict our search for the optimal solution to the following {\em compact} set
\[\small \mathcal{C}_F:=\hspace{-0.1cm}\lbrace (\lambda, U): U\in\mathcal{U}, \xi \geq \lambda \geq\epsilon , (\lambda^{-1}\ofd(U)+\hat{\Sigma}_m)\succeq \beta I \rbrace\]
for a certain $\beta>0$.\\

\begin{teor}
Problem \eqref{dual} admits a solution.
\end{teor}
\proof
Since $\mathcal{C}_F$ is closed and bounded and $\tilde{J}$ is continuous over $\mathcal{C}$ and hence over $\mathcal{C}_F$, by Weirstrass's Theorem the minimum exists.
\qed

\subsection{Uniqueness}

{\vale It is worth noting that $\tilde J$ is convex over $\mathcal{C}_F$,
however, as we will show, it is not strictly convex.} Therefore establishing the uniqueness of the minimum is not a trivial task.

To streamline the notation, let us define $\tilde U:= \ofd(U)$.
Then, the following Proposition characterizes the second variation of $\tilde{J}$ in  direction
$(\delta \lambda, \delta \tilde{U})$, i.e. $\delta^{2}\tilde{J} (\lambda,\tilde{U};\delta \lambda, \delta \tilde{U})$.
\begin{propo} \label{prop_hessian}
Let  $\tilde{u}=\text{vec}(\tilde{U})$, $\delta \tilde{u}:=\text{vec}(\delta \tilde{U})$, and $K:= (\lambda^{-1}\tilde{U}+\hat{\Sigma}_m)^{-1}\otimes (\lambda^{-1}\tilde{U}+\hat{\Sigma}_m)^{-1}$. Let also
\[ H:= \left[ \begin{array}{cc}
\lambda^{-3}\tilde{u}\tp K \tilde{u} & -\lambda^{-2}\tilde{u}\tp K  \\
-\lambda^{-2}K\tilde{u} & \lambda^{-1}K\\
\end{array} \right]\in \mathbb{R}^{(1+n^2) \times (1+n^2)}. \]
Then, we have
$$
\delta^{2}\tilde{J} (\lambda,\tilde{U};\delta \lambda, \delta \tilde{U})=[\delta \lambda\ \ \delta \tilde{u}\tp]
H\bmat{c}\delta \lambda\\ \delta \tilde{u}\emat.
$$
\end{propo}
{\vale The proof is an easy computation and it is left to the reader.}

\begin{corr}\label{direction}
The functional $\tilde{J}$ is convex and for any point $(\lambda_0,\tilde{U}_0)$ there is exactly one direction along which it is not strictly convex. This direction is
\beq
\label{direzionemaledetta}
(\delta\lambda_0,\delta\tilde{U}_0)=(h\lambda_0,h\tilde{U}_0),\; \; \hbox{with } h\neq 0.
\eeq
\end{corr}
\proof
Since in $\mathcal{C}_F$ we have that $K\in {\mathbf Q}_{n^2}$ is positive definite and $\lambda>0$, the Hessian matrix  $H$ has at least rank equal to $n^2$.
Hence there is at most one direction along which $\tilde{J}$ is not strictly convex.
By direct computation it is immediate to check that the second variation along the direction \eqref{direzionemaledetta}
is zero.
\qed
Next we show that if at a certain point the functional $\tilde{J}$ is constant  along an arbitrary direction then $\tilde{J}$ vanishes in this point. 
\begin{lemm}\label{corstconvex}
Let $(\lambda_0,\tilde{U}_0)$ be a given point in the feasible set $\mathcal{C}_F$.
If $w:=(\delta\lambda,\delta \tilde{U})\neq(0,0)$ is any direction along which $\tilde{J}(\lambda_0,\tilde{U}_0)$ is constant, that is if there exists $\varepsilon>0$ such that
$f(\alpha):=\tilde{J}(\lambda_0+\alpha \delta\lambda,\tilde{U}_0+\alpha \delta \tilde{U})$ is constant for any $\alpha$ such that $|\alpha|<\varepsilon$,  then $\tilde{J}(\lambda_0,\tilde{U}_0)=0$.
\end{lemm}
\proof
By assumption we have that  $f(\alpha)$ is constant in a neighbourhood of zero.
Hence the first derivative $f'(0)$, and thus the second derivative $f''(0)$, must vanish and hence the second variation
of $f(0)$ in direction $1$ is zero.
On the other hand this second variation is, by definition, the second variation $\delta^2\tilde{J}(\lambda_0,\tilde{U}_0,\delta\lambda, \delta \tilde{U})$.
Hence this second variation vanishes and,
by Corollary \ref{direction}, this implies that 
$$(\delta\lambda, \delta \tilde{U})=(h\lambda_0,h\tilde{U}_0),$$
for a certain real constant $h$.
Hence, for $|\alpha|$ sufficiently small, we have 
\beq\label{semif}
\tilde{J}(\lambda_0,\tilde{U}_0)=f(0)=f(\alpha)=\tilde{J}((1+\alpha h)\lambda_0,(1+\alpha h)\tilde{U}_0).
\eeq
By direct computation we get
$$
\tilde{J}((1+\alpha h)\lambda_0,(1+\alpha h)\tilde{U}_0)=(1+\alpha h)\tilde{J}(\lambda_0,\tilde{U}_0)
$$
which together with (\ref{semif}) yields the conclusion.
\qed

We are now ready to state our main result.
\begin{teor}
The dual problem \eqref{dual} admits a unique solution.
\end{teor}
\proof
By contradiction, assume that there exist two optimal solutions $(\lambda^{*}_1,\tilde{U}_1^{*})$ and $(\lambda^{*}_2,\tilde{U}_2^{*})$. By the convexity of the set $\mathcal{C}_F$, the whole segment ${\mathcal S}$ connecting $(\lambda^{*}_1,\tilde{U}_1^{*})$ to $(\lambda^{*}_2,\tilde{U}_2^{*})$ must belong to $\mathcal{C}_F$. 
Then, by the convexity of $\tilde{J}(\cdot, \cdot)$ all the points in  ${\mathcal S}$ are optimal solutions so that  $\tilde{J}(\cdot, \cdot)$ is constant in  ${\mathcal S}$.
In view of  Lemma \ref{corstconvex} this implies that $\tilde{J}(\cdot, \cdot)$ is zero in  ${\mathcal S}$ and this is a contradiction since {\vale it can be proved that} the optimal value of $\tilde{J}$ is negative.
\qed

\begin{corr}\label{optonbound}
Any optimal solution $(\lambda^{*},\tilde{U}^{*})$ minimizing $\tilde{J}$ over $\mathcal{C}_F$  lies on the boundary of $\mathcal{C}_F$.
\end{corr}
\proof
Let $(\lambda^{*},\tilde{U}^{*})$ be an optimal solution and, by contradiction, assume that
$(\lambda^{*},\tilde{U}^{*})$ does not belong to the boundary of $\mathcal{C}_F$. Then there exists $\varepsilon>0$ such that  $$((1+\varepsilon)\lambda^{*},(1+\varepsilon)\tilde{U}^{*})\in \mathcal{C}_F.$$
Now by direct computation 
\beq
\tilde{J}((1+\varepsilon)\lambda^{*},(1+\varepsilon)\tilde{U}^{*})=(1+\varepsilon)\tilde{J}(\lambda^{*},\tilde{U}^{*})<\tilde{J}(\lambda^{*},\tilde{U}^{*})
\eeq
where the last inequality follows from the fact that, as 
seen in the proof of Lemma \ref{lemma_l0}, the optimal value of $\tilde{J}$ is negative.
This a contradiction since $\tilde{J}(\lambda^{*},\tilde{U}^{*})$ is assumed to be a minimum. 
\qed
{\vale \begin{rem} 
The zero duality gap between the primal and the dual problem can be exploited to recover the solution of the primal problem \eqref{problem}, for more details see e.g. \cite{zorzi2016ar}, \cite{ciccone2017factor}.
\end{rem}}

\section{Conclusion}
In this paper the problem of robust latent-variables graphical model identification has been considered.
In particular, the \textit{Sparse plus Low Rank} decomposition problem has been reformulated for the case in which only the sample covariance is available and the difference between the sample covariance and the actual one is non-negligible. 
The extension of the presented work to the dynamical case will be subject for future investigation.



%


\end{document}